\documentclass{article}
\usepackage{latexsym, fullpage}

\usepackage{mathrsfs}
\usepackage[english]{babel}
\usepackage[utf8]{inputenc}
\usepackage{graphicx, color}
\usepackage[colorinlistoftodos]{todonotes}
\usepackage{fullpage}
\usepackage{graphics, color}
\usepackage{psfrag}
\usepackage{amsmath}
\usepackage{enumerate}
\usepackage{amsthm}
\usepackage{hyperref}
\usepackage{calrsfs}
\usepackage{tikz}
\usepackage{rotating}
\usepackage{amssymb}
\usepackage{amsmath}

\usepackage{authblk}
\usetikzlibrary{positioning,chains,fit,shapes,calc}

\newtheorem{theorem}{Theorem}[section]

\newtheorem{lemma}[theorem]{Lemma}

\newtheorem{prop}[theorem]{Proposition}
\newtheorem{claim}{Claim}
\theoremstyle{remark}

\title{The spectral radius of $1$-planar graphs without complete subgraphs}

\author{{Weilun Xu\thanks{Email: WeilunXu94@163.com}}, {An Chang\thanks{Email: anchang@fzu.edu.cn}}
\\
\small{Center for Discrete Mathematics and Theoretical Computer Science\\ Fuzhou University\\ Fuzhou
Fujian, China}
 }

\date{ }

\begin{document}
\maketitle
\begin{center}
{\bf ABSTRACT}
\end{center}
A 1-planar graph refers to a graph that can be drawn on the plane such that each edge
has at most one crossing. In this paper, focusing on the spectral Tur\'{a}n-type problems of $1$-planar graphs, we determine completely the unique spectral extremal graph among all $K_3$-free or $K_4$-free $1$-planar graphs, and provide a characterization of the spectral extremal graphs for $K_5$-free $1$-planar graphs, confining the candidates to a specific, small family.


\section{Introduction}

All graphs considered in this paper are finite undirected graphs without loops and multiple edges.
Let $G$ be a graph with vertex set $\{1,2,\ldots,n\}$, and the \emph{adjacency matrix} of $G$ is the $n\times n$ matrix $A_G=(a_{ij})$, where $a_{ij}$ is 1 if $ij$ is an edge of $G$, and $0$ otherwise.
The $spectral$ $radius$ of $G$ is the largest eigenvalue of its adjacency matrix $A_G$. One of the most extensively studied problems in spectral graph theory concerns the relation between the spectral radius and various properties of a graph, of which Brualdi and Solheid \cite{Brualdi} presented the following well-known problem in 1986:

{\textbf{Problem 1}. (Brualdi–Solheid problem) Given a set $\mathcal{G}$ of graphs, find a tight upper bound for the spectral radius in $\mathcal{G}$ and characterize the extremal graphs meeting it.

Numerous classes of graphs have been well investigated since then, leading to many excellent works in the literature, such as graphs with cut vertices \cite{Berman}, given chromatic index \cite{Feng}, given diameter \cite{Hansen}, and so on.

We use $P_k$, $C_{k}$ to denote a path with $k$ vertices, a cycle with $k$ vertices. Let $G_1$ and $G_2$ be two graphs with disjoint vertex sets. The \emph{join} of $G_1$ and $G_2$ is denoted by $G_1+G_2$, that is, $V(G_1+ G_2)=V(G_1)\cup V(G_2)$ and $E(G_1+ G_2)=E(G_1)\cup E(G_2)\cup \{xy:x\in V(G_1)~\text{and}~y\in V(G_2)\}$. The \emph{Cartesian product} $G \mathbin{\square} H$ of graphs $G$ and $H$ has vertex set $V(G) \times V(H)$, with $(u,v) \sim (u',v')$ if and only if either $u = u'$ and $vv' \in E(H)$, or $v = v'$ and $uu' \in E(G)$.
For a given graph class $\mathcal{H}$, a graph $G$ is said to be $\mathcal{H}$-free if $G$ does not contain any graph from $\mathcal{H}$ as a subgraph. In particular, when $\mathcal{H}=\{H\}$, we simplify the notation to $H$-free. If the class of graphs $\mathcal{G}$ in Problem 1 is restricted to an $\mathcal{H}$-free graph for a given graph class $\mathcal{H}$, it will result in the following so-called spectral Tur\'{a}n-type problem:

{\textbf{Problem 2}. (Spectral Tur\'{a}n-type problem) What is the maximum spectral radius of an $\mathcal{H}$-free graph
$G$ with $n$ vertices for a given graph class $\mathcal{H}$?

Analogous to the classical Tur\'{a}n-type problems, the case when $\mathcal{H}$ contains only a complete graph is one of the most fundamental questions in the study of the spectral Tur\'{a}n-type problem. This problem was solved by Nikiforov, a pioneering figure in the study of spectral Tur\'{a}n problems. In \cite{Nikiforov Kr}, he proved that among all $n$-vertex $K_r$-free graphs, the complete $(r-1)$-partite graph with parts of equal or nearly equal size attains the maximum spectral radius. Since then, a large number of remarkable results on spectral Tur\'{a}n-type problems have been achieved. For more on this issue, we refer the reader to the nice survey \cite{Li}.

We emphasize the noteworthy work of Tait and Tobin \cite{Tait}, who in 2017 examined the Brualdi–Solheid problem for the class of planar graphs on
$n$ vertices. They demonstrated that, for $n$ large enough, the unique graph attaining the maximum spectral radius among all $n$-vertex planar graphs is $K_2+P_{n-2}$. The approach employed by Tait and Tobin \cite{Tait} in their proof has gained recognition as the \emph{second characteristic equation method}. This method offers important insights for the study of both Problem 1 and Problem 2.

A natural question arising from the results of Tait and Tobin is the following planar spectral Tur\'{a}n-type problem.

{\textbf{Problem 3}. (Planar spectral Tur\'{a}n-type problem) Let $\mathcal{H}$ be a given class of planar graphs. Among all $n$-vertex $\mathcal{H}$-free planar graphs, which one attains the maximum spectral radius?

Let $SPEX_{\mathcal{P}}(n,\mathcal{H})$ be the set of graphs that maximize the spectral among all $n$-vertex $\mathcal{H}$-free planar graphs.  To our knowledge, the first result on this problem was obtained by Zhai and Liu \cite{Zhai}, who determined $SPEX_{\mathcal{P}}(n,\mathcal{H})$ when $\mathcal{H}$ consists of $k$ edge-disjoint cycles. Subsequently, the $SPEX_{\mathcal{P}}(n, H)$ was determined for many specific graphs $H$, including $t$ vertex-disjoint cycles \cite{Fang}, two cycles sharing one common edge \cite{ZhangW}, the Theta graph \cite{ZhangW}, the wheel graph \cite{WangH}, the friendship graph \cite{WangH}, the book graph \cite{Xu}, among others.

Similar to Problem 3, the question that naturally arises from the conclusions of Tait and Tobin \cite{Tait} is: In which classes of graphs more general than planar graphs can the graph achieving the maximum spectral radius be determined? In light of this question, one can investigate the   classes of graphs that resemble planar graphs the most, i.e., $1$-planar graphs.

A $drawing$ of a graph $G=(V,E)$ is defined as a mapping $D:G\to \mathbb{R}^2$ that uniquely assigns a distinct point in the plane to each vertex in $V$ and a continuous arc connecting $D(u)$ and $D(v)$ to each edge $uv \in E$. Furthermore, a drawing is termed a $1$-$planar$ $drawing$ if each of its edges $D(uv)$ is crossed at most once. Consequently, a graph $G$ is designated as a $1$-$planar$ $graph$ if it admits a $1$-planar drawing. A $1$-planar graph $G$ together with a $1$-planar drawing is a $1$-$plane$ $graph$.
The concept of $1$-planar graphs was originally introduced by Ringel\cite{Ringel}, and since then many properties of 1-planar graphs have been studied \cite{Kobourov}.

Although the class of $1$-planar graphs looks 'similar' to that of planar graphs, there are profound differences between the two classes. A distinct difference lies in the fact that any planar graph is $\{K_5,K_{3,3}\}$-minor-free, which is a consequence of the famous Kuratowski's theorem. However, for any graph $H$, we can construct a $1$-planar graph that contains an $H$-minor by subdividing the edges of $H$. Therefore, the class of $1$-planar graphs cannot be characterized solely by forbidden finite minors. Moreover, by Euler's formula, an $n$-vertex edge-maximal planar graph has exactly $3n-6$ edges. Consequently, the class of edge-maximal planar graphs is identical to the class of edge-maximum planar graphs. However, there exists an $n$-vertex edge-maximal $1$-planar graph that has at most $\frac{45n}{17}+O(1)$ edges \cite{Brandenburg}, significantly fewer than the maximum number of edges possible in an $n$-vertex $1$-planar graph, which is $4n-8$ \cite{Pach}.

In 2024, Zhang, Wang and Wang \cite{ZhangWW} considered the spectral extremal problem on $1$-planar graphs. They proved that for sufficiently large $n$, the unique $n$-vertex $1$-planar graph attaining the maximum spectral radius is $K_2 + P_{n-2}^{2+}$, where $P_{n}^{2+}$ denotes the graph obtained from a path $v_1, v_2, v_3, \dots, v_n$ by adding the edge $v_1v_n$ and edges $v_iv_{i+2}$ for all $1 \leq i \leq n-2$. (See Figure \ref{cpn}).
\begin{figure}[h!]
    \centering
    \includegraphics[width=0.4\textwidth]{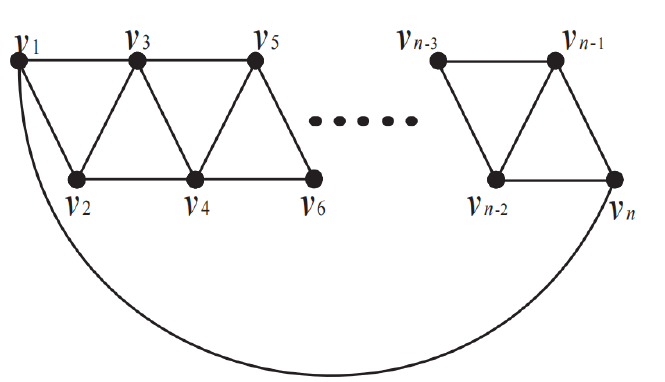}
    \caption{The graph $P_{n}^{2+}$.}
    \label{cpn}
\end{figure}

Just as the conclusion of Tait and Tobin \cite{Tait} gives rise to the planar spectral Tur\'{a}n problem, so does the result of Zhang, Wang, and Wang lead to the following $1$-planar spectral Tur\'{a}n problem.

{\textbf{Problem 4}. ($1$-planar spectral Tur\'{a}n-type problem) Let $\mathcal{H}$ be a given class of $1$-planar graphs. Among all $n$-vertex $\mathcal{H}$-free $1$-planar graphs, which attains the maximum spectral radius?

 Denote $SPEX_{\mathcal{P}_1}(n,H)$ the set of graphs that achieve the maximum spectral radius among all $n$-vertex $H$-free $1$-planar graphs. In this paper, we consider the fundamental case of Problem 4, i.e., $H$ is a complete graph. Since $K_2+P_{n-2}^+$ contains no $K_t$ when $t\geq 6$, it follows that for $t\geq 6$, $SPEX_{\mathcal{P}_1}(n,K_t) = \{K_2+P_{n-2}^+\}$ from the main result mentioned above by Zhang, Wang and Wang in \cite{ZhangWW}. Therefore, it remains to consider the case $H \in \{K_3,K_4, K_5\}$. When $H=K_3$, we give the following theorem as an immediate consequence of Lemma 2.7 in Section \ref{Pre}.
 \begin{theorem}
 Let $n$ be a sufficiently large integer. Then $SPEX_{\mathcal{P}_1}(n,K_3) = \{K_{2,n-2}\}$.
 \end{theorem}

 Define the graph $C_{\frac{n}{2}}^{\square}$ as follows (see Figure \ref{Cn} for example):\\
(1) If $n$ is even, then $C_{\frac{n}{2}}^{\square}$ is the Cartesian product of $C_{\frac{n}{2}}$ and $K_2$.\\
(2) If $n$ is odd, then $C_{\frac{n}{2}}^{\square}$ is obtained from $C_{\frac{n}{2}}^{\square}$ in (1) by splitting any one vertex, say $v$, of $C_{\frac{n}{2}}^{\square}$, i.e., to replace $v$ by two nonadjacent vertices of degree two, $v'$ and $v''$, so that two edges incident
to $v$ in the same $C_4$ are incident to either $v'$ or $v''$, the other end of these two edges remaining unchanged.
\begin{figure}[h!]
    \centering
    \includegraphics[width=0.5\textwidth]{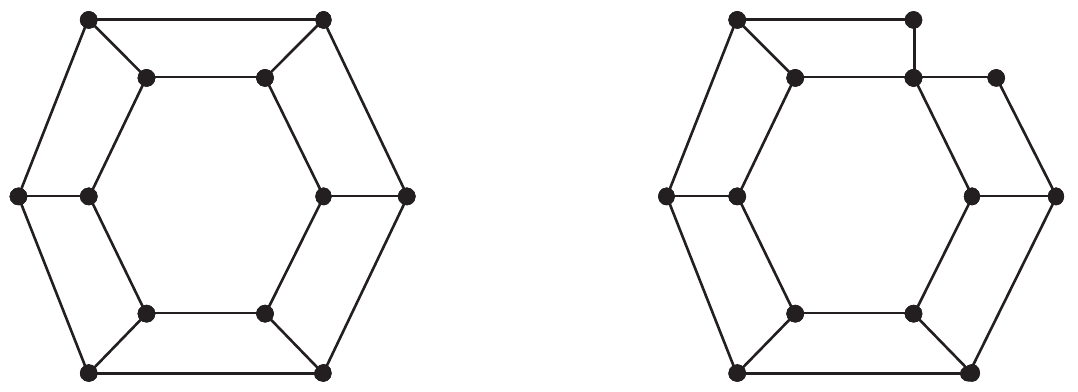}
    \caption{The graphs $C_{6}^{\square}$ and $C_{\frac{13}{2}}^{\square}$.}
    \label{Cn}
\end{figure}

 When $H=K_4$   We obtain the following result:
\begin{theorem}\label{Th-K4}
Let $n$ be a sufficiently large integer. Then $SPEX_{\mathcal{P}_1}(n,K_4) = \{2K_1+C_{\frac{n-2}{2}}^{\square}\}$.
\end{theorem}

Let $\mathcal{P}_n^2$ be the set of graphs obtained from $P_n^{2+}$ by deleting exactly  $\left\lceil\frac{n}{2}\right\rceil$ edges so that all triangles are destroyed.  Obviously, $\mathcal{P}_n^2$ contains exactly one graph $QP_n$ when $n$ is even, which is obtained from $P_n^{2+}$ by removing the edges $v_{2i}v_{2i+1}$ for all $1\leq i \leq \frac{n-2}{2}$. When $n$ is odd, the situation becomes more involved, and a precise characterization will be presented in Section \ref{K5}. Moreover, let $C_{n}^2$ be the graph obtained from a $C_n=v_1v_2\dots v_nv_1$ by adding edges $v_iv_{(i+2\mod n)}$ for all $1\leq i \leq n$. In particular, let $C_n^{2-}$ be the graph obtained from $C_{n}^2$ by removing $v_{n}v_2$. We now present our result on $H=K_5$ as follows.

\begin{theorem}\label{Th-K5}
Let $n$ be a sufficiently large integer. If $n$ is even, then $SPEX_{\mathcal{P}_1}(n,K_5) \subseteq \{2K_1+C_{n-2}^{2},K_2+QP_{n-2}\}$. If $n$ is odd, then $SPEX_{\mathcal{P}_1}(n,K_5)\subseteq \{2K_1+C_{n-2}^{2-}\}\cup \{K_2+G:G\in \mathcal{P}_{n-2}^2 \}$.
\end{theorem}

The rest of this paper is organized as follows. In Section \ref{Pre}, we introduce the notation and terminologies which will be used frequently. In particular, a key lemma will be proven in this section. In Section \ref{K4}, we will  present the proof of Theorem \ref{Th-K4}. In Section \ref{K5}, we will give the proof of Theorem \ref{Th-K5}.

\section{Preliminaries}\label{Pre}
Throughout the subsequent discussion, when we refer to a graph $G=(V,E)$ with $n$ vertices, we always assume that its vertex set $V$ is $\{1,2,3,\dots, n\}$. For a subset $A$ of $V$, we write $G[A]$ for the subgraph of $G$ induced by $A$, that is, the graph obtained from $G$ by removing vertices in $V\setminus A$. Let $e\in E$ and $e'\in V^2\setminus E$. We use $G-e$ and $G+e'$ to denote the graphs obtained from $G$ by removing $e$ and adding $e'$, respectively. Let $A_G$  be the adjacency matrix of $G$, and
the spectral radius of $G$ is denoted by $\lambda(G)$. For an $n$-dimensional vector $\mathbf{v} = (\mathbf{v}_1, \mathbf{v}_2, \ldots, \mathbf{v}_n)^{\intercal}\in \mathbb{R}^n$, we define the maximum norm ($l_{\infty}$-norm) for this vector as $\|\mathbf{v}\|_{\infty} = \max_{1 \leq i \leq n} |\mathbf{v}_i|$.

For a connected graph $G$ on $n$ vertices and a positive constant $c$, the famous Perron-Frobenius Theorem in the theory of nonnegative matrices guaranties that there exists a unique eigenvector in $\{\mathbf{v} \in \mathbb{R}^{n} : \|\mathbf{v}\|_{\infty} = c\}$ corresponding to the spectral radius $\lambda(G)$,  and it has strictly positive entries. We refer to this eigenvector as \emph{the Perron vector} of $G$.

Clearly, $A_G$ is a symmetric matrix for any graph $G$. Thus, the spectral radius of a graph $G$ can be obtained from the well-known Courant-Fischer Theorem, which states that
\[
\lambda(G) = \max_{\mathbf{v} \in \mathbb{R}^n} \frac{\mathbf{v}^{\intercal} A_G \mathbf{v}}{\sum_{i=1}^n\mathbf{v}_i^2}.
\]
In what follows, we will employ the following lemma to compare the spectral radius of two graphs, which is a straightforward corollary of the Courant-Fischer Theorem.
\begin{lemma}\label{lem-3}
Let $G$ and $G'$ be two $n$-vertex graphs, and let $\mathbf{x}$ be the Perron vector of $G$. If
$$\mathbf{x}^{\intercal} A_{G'} \mathbf{x}>\mathbf{x}^{\intercal} A_G \mathbf{x},$$
then $\lambda(G')>\lambda(G).$
\end{lemma}

Next, we present several structural properties of $1$-planar graphs. In particular, this paper distinguishes between $1$-planar graphs and $1$-plane graphs. When referring to the specific $1$-plane graph corresponding to a $1$-planar graph $G$, we denote it by $D(G)$, where $D$ represents the $1$-planar drawing associated with $G$.
\begin{prop}\cite{Fabrici}\label{pro-4}
Every $n$-vertex $1$-planar graph has at most $4n-8$ edges.
\end{prop}

\begin{prop}\cite{Fabrici}\label{pro-1}
Every $1$-planar graph is $7$-degenerate.
\end{prop}
\begin{prop}\cite{Huang}\label{pro-2}
Every $1$-planar graph is $K_{3,7}$-free.
\end{prop}
\begin{prop}\label{pro-3}
Let $D(G)$ be a $1$-plane graph, and let $u_1v_1$, $u_2v_2$ be a pair of crossing edges. Then for any pair of vertices $s\in\{u_1,v_1\}$ and $t\in\{u_2,v_2\}$, there exists a connected region $\mathcal{R}_{st}$ such that $\mathcal{R}_{st}\cap D(G)=\{s,t\}$.
\end{prop}

\begin{lemma}\cite{ZhangWW}\label{lem-1}
Let $\mathcal{G}$ be a class of $1$-planar graphs and $G$ be the graph that maximizes the spectral radius among all $n$-vertex graphs in $\mathcal{G}$. Suppose that $\mathbf{x}$ is the Perron vector of $G$ with $\|\mathbf{x}\|_{\infty}=1 $. If $\sqrt{2n-4}\leq\lambda(G)$ and for any $v\in V(G)$, $\mathbf{x}_v \geq \frac{1}{\lambda(G)}$ for sufficiently large $n$, then for a constant $\epsilon \leq \frac{1}{21000}$ independent of $n$, the following hold:\\
(1) Let $x$ be a vertex such that $\mathbf{x}_x=1$. Then $d_x\geq (1-24\epsilon)n$. Moreover, there exists a vertex $w\neq x$ such that $\mathbf{x}_w\geq 1-47\epsilon$ and $d_w\geq (1-116\epsilon)n$.\\
(2) For any $v\in V(G)\setminus\{x,w\}$, $\mathbf{x}_v\leq \frac{1}{20}$.

\end{lemma}

\begin{lemma}\label{lem-2}
Suppose $H\in\{K_3,K_4,K_5\}$ and $G\in SPEX_{\mathcal{P}_1}(n,H)$. If $n$ is sufficiently large, then $G$ contains a complete bipartite graph $K_{2,n-2}$ as a spanning subgraph.
\end{lemma}
\proof
Let $\mathbf{x}$ be the Perron vector of $G$ with $\|\mathbf{x}\|_{\infty} = 1$ and $x$ be a vertex such that $\mathbf{x}_x=1$.

We first show that $\sqrt{2n-4}\leq\lambda(G)$. Clearly, $K_{2,n-2}$ is an $H$-free $1$-planar graph and thus $\lambda(G)\geq \lambda(K_{2,n-2})=\sqrt{2n-4}$. Let $v$ be a vertex of $G$. If $\mathbf{x}_v<\frac{1}{\lambda(G)}$, then the graph $G'$ obtained from $G$ by removing all edges incident to $v$ and adding the edge $xv$ remains $H$-free as $v$ has degree $1$ in $G'$. Moreover,
$$\lambda(G')\mathbf{x}^{\intercal}\mathbf{x}\geq\mathbf{x}^{\intercal} A_{G'}\mathbf{x}\geq \mathbf{x}^{\intercal} A_{G}\mathbf{x} + 2\mathbf{x}_v-2\mathbf{x}_v\sum_{u\sim v}\mathbf{x}_u =\lambda(G)\mathbf{x}^{\intercal}\mathbf{x} + 2\mathbf{x}_v(1-\lambda(G)\mathbf{x}_v)>\lambda(G)\mathbf{x}^{\intercal}\mathbf{x},$$
which contradicts the maximality of $\lambda(G)$.

Now by Lemma \ref{lem-1} (1), $d_x\geq (1-24\epsilon)n$, and  there exists a vertex $w\neq x$ such that $\mathbf{x}_w\geq 1-47\epsilon$ and $d_w\geq (1-116\epsilon)n$. Thus, $|N(x)\cap N(w)|\geq (1-24\epsilon + 1 - 116\epsilon)n-n\geq (1-140\epsilon)n$, where $\epsilon$ is a constant defined as in Lemma \ref{lem-1}. Let $A= N(x)\cap N(w)$ and $B=V(G)\setminus (A\cup\{x,w\})$. Then $|B|\leq 140\epsilon$ and $\mathbf{x}_v\leq \frac{1}{20}$ for any $v\in A\cup B$.

Next, we will show that $B=\emptyset$. Suppose, to the contrary, that $B\neq \emptyset$. Then by Proposition \ref{pro-1}, there exists a sequence $v_1v_2,\dots v_{|B|}$ composed of vertices in $B$ such that $v_{i}$ has at most $7$ neighbors in $\{v_{i+1},v_{i+2},\dots, v_{|B|}\}$ for all $1\leq i \leq |B-1|$. Now for $1\leq i \leq |B|-1$, define $G_1=G$ and $G_{i+1}$ is obtained form $G_i$ by removing all the edges incident to $v_{i}$ and adding the edges $xv_{i}$, $wv_{i}$. Note that for any $1\leq i \leq |B-1|$, $v_{i}$ has at most $15$ neighbors in $G_i$: at most $1$ in $\{x, w\}$, at most $7$ in $B$ and at most $7$ in $A$, which follows from the definition of $B$ and Proposition \ref{pro-2}, respectively. Thus, by Lemma \ref{lem-1} (2), we have
$$\mathbf{x}^{\intercal} A_{G_{i+1}}\mathbf{x}-\mathbf{x}^{\intercal} A_{G_{i}}\mathbf{x}\geq 2((1+1-47\epsilon)\mathbf{x}_{v_{i}}-(1+\frac{7}{10})\mathbf{x}_{v_{i}})>\frac{2}{5}\mathbf{x}_{v_{i}}>\frac{2}{5\lambda(G)}.$$

Let $G'=G_{|B|}$. If $B\neq \emptyset$, then we have
\begin{equation}\label{eq-1}
\lambda(G')-\lambda(G)>\frac{2}{5\lambda(G)}.
\end{equation}
Note that any vertex in $B$ has degree $2$ in $G'$. Thus, if $H\in\{K_4,K_5\}$, then $G'$ is $H$-free. If $H=K_3$, then $x\nsim w$ as $A\neq \emptyset$. Thus, $G'$ is still $H$-free. Consequently, $G'$ is not a $1$-planar graph.

By the definition of $G'$, $G'[A\cup\{x,w\}]=G[A\cup\{x,w\}]$. Thus, $G'[A\cup\{x,w\}]$ is a $1$-planar graph as $G$ is a $1$-planar graph.
Let $D(G'[A\cup\{x,w\}])$ be a $1$-plane graph. If there exist two vertices $u,v\in A$ such that $xu$ crosses $wv$, then by Proposition \ref{pro-3}, there exists a region $\mathcal{R}_{xw}$ such that $\mathcal{R}_{xw}\cap D(G'[A\cup\{x,w\}])=\{x,w\}$. However, in this case, we can extend the $1$-planar drawing of $G'[A\cup\{x,w\}]$ to $G'$ by drawing all vertices of $B$  within $\mathcal{R}_{xw}$ such that $D(G'[B\cup\{x,w\}])$ is a plane graph, which is a contradiction.

Now choose a $1$-planar drawing of $G'[A\cup\{x,w\}]$ and denote the corresponding $1$-plane graph by $D(G'[A\cup\{x,w\}])$. From the above discussion, the graph induced by the edge set $\{e:e\cap\{x,w\}\neq\emptyset\}$, denoted by $D(G_1)$, is a plane graph. In $D(G_1)$, let the vertices in $A$ be labeled $v_1, v_2, \dots, v_{|A|}$ in counterclockwise order around $x$. For each $i$ (with indices modulo $|A|$), let $f_i$ be the face incident to the vertices $v_{i}$, $v_{i+1}$, $x$, and $w$. In particular, if $xw\in E(G)$, then let $f_j$ be the unique face that contains edge $xw$ in $D(G'[A\cup\{x,w\}])$. Obviously, for any $1\leq i\leq |A|$, $N(v_i)\cap A\subseteq\{v_{i-2},v_{i-1},v_{i+1},v_{i+2}\}$.
Hence, for $i\neq j$, there are at most three possible edges,  $v_{i-1}v_{i+1}$, $v_{i}v_{i+1}$, and $v_{i}v_{i+2}$ in $D(G'[A\cup\{x,w\}])$ that intersect the interior of $f_i$, and therefore $G'-\{v_{i-1}v_{i+1},v_{i}v_{i+1},v_{i}v_{i+2}\}$ is a $1$-planar graph because we can draw  all vertices of $B$  within $f_i$ such that $D(G'[B\cup\{x,w\}])$ is a plane graph.

It remains to show that there exists a suitable $i$. An ordering $4$-tuple $(i-1,i,i+1,i+2)$ is called $good$ if $\max\{\mathbf{x}_{v_{i-1}},\mathbf{x}_{v_{i}},\mathbf{x}_{v_{i+1}},\mathbf{x}_{v_{i+2}}\}\leq\frac{80}{\lambda(G)}$ and $bad$ otherwise. Define $L=\{v\in A: \mathbf{x}_v>\frac{80}{\lambda(G)}\}$. Then we have
$$2(4n-8)\geq \sum_{v\in V(G)}d_v\geq\sum_{v\in L}d_v\geq \sum_{v\in L}\lambda(G)\mathbf{x}_v\geq 80|L|,$$
where the first inequality follows from  Proposition \ref{pro-4}, and the third inequality follows from the characteristic equation. Thus, $|L|\leq \frac{n}{10}$. Note that each vertex in $L$ gives rise to exactly $4$ bad ordering $4$-tuples. Recall that $|A|\geq (1-140\epsilon)n$, and thus there are at least $(1-140\epsilon)n$ ordering $4$-tuples. Therefore, there exist at least $(1-140\epsilon)n-\frac{2}{5}n>2$ good ordering $4$-tuples. Pick $i\neq j$ such that $(i-1,i,i+1,i+2)$ is good, and remove edge set $\{v_{i-1}v_{i+1},v_{i}v_{i+1},v_{i}v_{i+2}\}$ from $G'$. Denote the result graph by $G''$. Clearly, $G''$ is an $H$-free $1$-planar graph. Moreover,
\begin{equation}\label{eq-2}
\mathbf{x}^{\intercal} A_{G'}\mathbf{x}-\mathbf{x}^{\intercal} A_{G''}\mathbf{x}\leq 6\times\frac{6400}{\lambda(G)^2}.
\end{equation}
Recall that $n$ is sufficiently large and $\sqrt{2n-4}\leq\lambda(G)$.  Combining inequalities (\ref{eq-1}) and (\ref{eq-2}), we have
$$\mathbf{x}^{\intercal} A_{G''}\mathbf{x}-\mathbf{x}^{\intercal} A_{G}\mathbf{x}\geq \frac{2}{5\lambda(G)}-\frac{38400}{\lambda(G)^2}>0,$$
which contradicts the maximality of $\lambda(G)$.
\endproof
A direct consequence of Lemma \ref{lem-2} is that $SPEX_{\mathcal{P}_1}(n,K_3)=\{K_{2,n-2}\}$ for sufficiently large $n$, that is the conclusion of our Theorem 1.1.

\section{Proof of Theorem \ref{Th-K4}}\label{K4}
To begin our proof, we will need the following lemma.
\begin{lemma}\label{lemplane}\cite{Xu}
Let $n$ be a sufficiently large integer. Among all $n$-vertex $K_4$-free planar graphs, $2K_1+C_{n-2}$ is the unique graph that reaches the maximum spectral radius.    
\end{lemma}

Let $G\in SPEX_{\mathcal{P}_1}(n,K_4)$ and $\mathbf{x}$ be  the Perron vector of $G$ with $\|\mathbf{x}\|_{\infty} =1$. Then by Lemma \ref{lem-2}, $K_{2,n-2}\subseteq G$. Denote $x$ and $w$ the two vertices in $G$ such that $d_G(x)\geq d_G(w)\geq n-2$, and denote $A=V(G)\setminus\{x,w\}$. Then, by Lemma \ref{lem-1}, we have $\mathbf{x}_x=\mathbf{x}_w=1$.

\begin{claim}
$xw\notin E(G)$.
\end{claim}
\proof
If $xw\in E(G)$, then $G[A]$ is an empty graph, as otherwise $G$ contains a $K_4$. Thus, $G=K_2+I_{n-2}$. Note that $G$ is a $K_4$-free planar graph. 
However, $\lambda(G) < \lambda(2K_1 + C_{n-2})$, which contradicts Lemma \ref{lemplane}.
\begin{claim}\label{cla-3}
For any $v\in A$, $\frac{2}{\lambda(G)}\leq\mathbf{x}_v\leq\frac{2}{\lambda(G)-7}.$
\end{claim}
\proof
Let $u$ be the vertex such that $\mathbf{x}_u=\max_{v\in A}\{\mathbf{x}_v\}$. By Proposition \ref{pro-2}, $u$ has at most $7$ neighbors in $A$. Thus, we have
$$\lambda(G)\mathbf{x}_u\leq 2+7\mathbf{x}_u,$$
which equals  $\mathbf{x}_u\leq\frac{2}{\lambda(G)-7}$.
On the other hand, for any $v\in A$, since $xv\in E(G)$ and $wv\in E(G)$,
$$\lambda(G)\mathbf{x}_v\geq 2.$$
Thus $\frac{2}{\lambda(G)}\leq\mathbf{x}_v\leq\frac{2}{\lambda(G)-7}$ as desired.
\begin{claim}\label{cla-2}
Let $G$ be an $n$-vertex $K_4$-free $1$-planar graph. If $H=K_{2,n-2}\subseteq G$, then there exists a $1$-planar drawing $D$ of $G$ such that $D(H)$ is a plane graph.
\end{claim}
\proof Let $D$ be a $1$-planar drawing of $G$ that minimizes the crossing number of $D(H)$. Denote the vertices of $H$ that have degree $n-2$ by $x$ and $w$, respectively. If $D(H)$ is a plane graph, then we are done. Now suppose that there exist $u,v\in V(G)$ such that $xu$ crosses $wv$.
In $D(H)$, let the vertices in $A$ be labeled $v_1,v_2,\dots,v_{n-2}$ in counterclockwise order around $x$. Moreover, by Proposition \ref{pro-3},
we can suppose that $xv_i$ crosses $wv_{j}$ only if $|i-j|=1$. Assume that $u=v_i$ and $v=v_{j}$. Then we have $|i-j|=1$. Suppose without loss of generality that  $v_1=u$, $v_2=v$ and the crossing point of  $xu$ and $wv$ is next to $v_{n-2}$. Then $N_G(u)\cup N_G(v)\subseteq \{x,w,u,v,v_3,v_4\}$.\\
{\bf Case 1.} $xv_3$ crosses $wv_4$ or $xv_4$ crosses $wv_3$.

By symmetry, we may assume that $xv_3$ crosses $wv_4$. If the crossing point of  $xv_3$ and $wv_4$ is in the region bounded by $x,v_4,v_3,w$, then there is no edge between $\{v_3,v_4\}$ and $\{u,v\}$, and $N_G(u)\cup N_G(v)\subseteq \{x,w,u,v\}$. We can modify $D$ to $D'$ such that  $D'(H[\{x,w,u,v\}])$ is a plane graph, which contradicts that $D$ minimizes the crossing number of $D(H)$. Thus, the crossing point of  $xv_3$ and $wv_4$ is not in the region bounded by $x,v_4,v_3,w$, and therefore $N_G(v_3)\cup N_G(v_4)\subseteq \{x,w,u,v\}$. We can modify $D$ to $D'$ such that  $D'(H[\{x,w,v_3,v_4\}])$ is a plane graph (see Figure \ref{rd1}), which contradicts the fact that $D$ minimizes the crossing number of $D(H)$.\\
\begin{figure}[h!]
    \centering
    \includegraphics[width=1\textwidth]{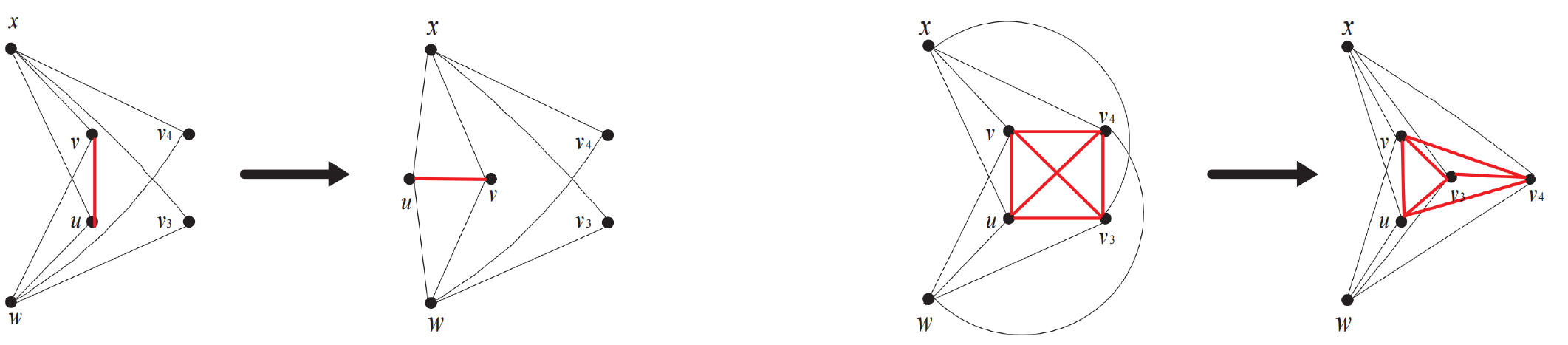}
    \caption{The redrawing of Case 1, where the red edges denote potential edges.}
    \label{rd1}
\end{figure}

\noindent{\bf Case 2.} $D(G[\{x,w,v_3,v_4\}]$) is a plane graph, and $v_4\in N_G(u)\cap N_G(v)$ .

In this case, it is clear that $v_3$ does not have a neighbor in $A \setminus \{u, v, v_4\}$. If $v_3v_4\in E(G)$, then we have $uv_3\notin E(G)$ and $wv_3\notin E(G)$. Let $G'=G-v_3v_4+v_3u+v_3v$. Then by Claim \ref{cla-3},
$$\mathbf{x}^{{\intercal}}A_{G'}\mathbf{x}-\mathbf{x}^{{\intercal}}A_{G}\mathbf{x}\geq 2\left(\frac{8}{(\lambda(G)-7)^2}-\frac{4}{\lambda(G)^2}\right)>0.$$
Moreover, it is easy to check that $G'$ is still $K_4$-free and $1$-planar, a contradiction. Thus, we have $v_3v_4\notin E(G)$.
And we can modify $D$ to $D'$ in such a way that $u$ and $v$ are located in the region bounded by $x,v_3,w,v_4$ and $D'(H[\{x,w,u,v\}])$ is a plane graph. (See Figure \ref{rd2}.) Obviously,  $D'(H[\{x,w,v_3,u,v,v_4\}])$ is a plane graph, which contradicts that $D$ minimizes the crossing number of $D(H)$.\\
\begin{figure}[h!]
    \centering
    \includegraphics[width=0.6\textwidth]{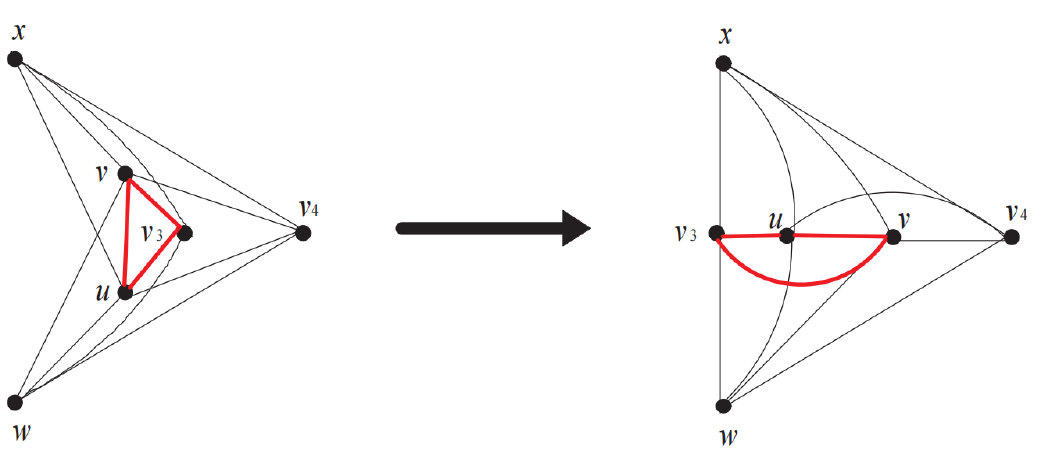}
    \caption{The redrawing of Case 2, where the red edges denote potential edges.}
    \label{rd2}
\end{figure}
{\bf Case 3.}$D(G[\{x,w,v_3,v_4\}]$) is a plane graph, and $|\{u,v\}\cap N_G(v_4)|\leq1$.

In this case, at most one of $u$ and $v$ is adjacent to $v_4$. Suppose without loss of generality that $uv_4\notin E(G)$. All possible edges of $G[\{u,v,v_3,v_4\}]$ are in the set $\{uv,uv_3,vv_3,v_3v_4, vv_4\}$. Then we can modify $D$ to $D'$ in such a way that $v$ is located in the region bounded by $x,v_3,w,u$ and $D'(H[\{x,w,u,v\}])$ is a plane graph. (See Figure \ref{rd3}.) Obviously,  $D'(H[\{x,w,u,v,v_3,v_4\}])$ is a plane graph, which contradicts that $D$ minimizes the crossing number of $D(H)$.\\
\begin{figure}[h!]
    \centering
    \includegraphics[width=0.7\textwidth]{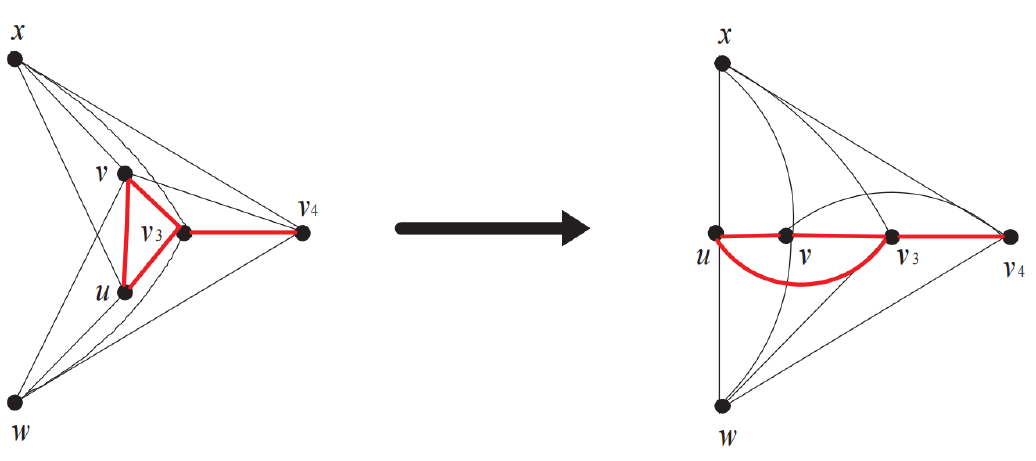}
    \caption{The redrawing of Case 3, where the red edges denote potential edges.}
    \label{rd3}
\end{figure}
According to cases 1--3, we complete the proof of claim.

By Claim \ref{cla-2}, let $D$ be the $1$-planar drawing of $G$ such that $D(H)$ is a plane graph where $H=K_{2,n-2}$. In $D(H)$, let the vertices in $A$ be labeled $v_1,v_2,\dots,v_{n-2}$ in counterclockwise order around $x$. Since $D(G)$ is a $1$-plane graph and $D(H)$ is a plane graph, we have $N_G(v_i)\subseteq \{x,w,v_{i-2},v_{i-1}, v_{i+1}, v_{i+2}\}~(i\mod{n-2})$. Thus,
$G[A]$ is a subgraph of $C_{n-2}^2$.  Moreover, since $G$ is $K_4$-free, $G[A]$ is obtained from $C_{n-2}^2$ by removing triangles. Let $\Delta_i$ be the triangle in $C_{n-2}^2$ formed by the three vertices $v_i$, $v_{i+1}$, and $v_{i+2}~(i\mod{n-2})$. Let $e_i=E(\Delta_i)\cap E(\Delta_{i-1})=v_{i}v_{i+1}~(i\mod{n-2})$ and $e_i'=E(\Delta_i)-e_i-e_{i+1}=v_{i}v_{i+2}~(i\mod{n-2})$. (See Figure \ref{fig2} for example.)
\begin{figure}[h!]
    \centering
    \includegraphics[width=0.6\textwidth]{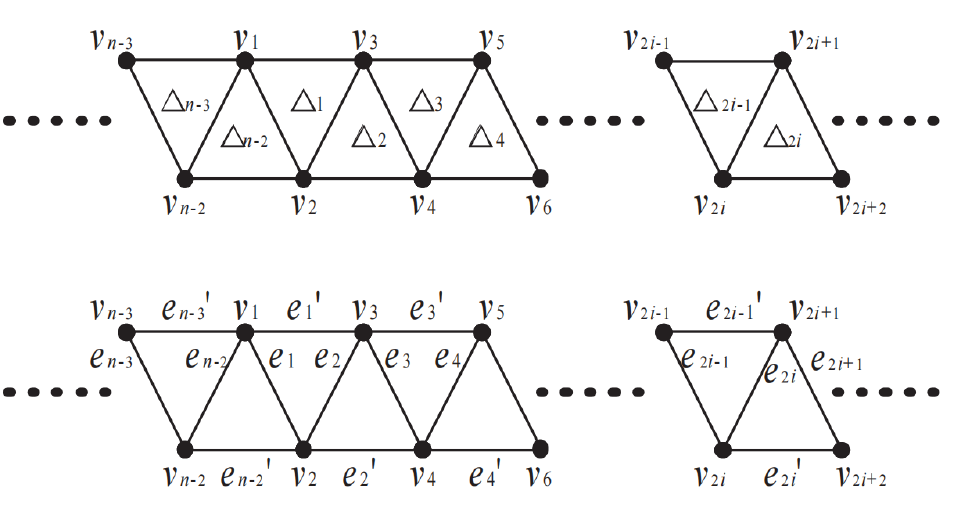}
    \caption{The supergraph of $G[A]$.}
    \label{fig2}
\end{figure}
Note that for any $\Delta_i$, there exists at least one edge in $E(\Delta_i)$ that is removed. Next, we consider two cases based on the parity of $n$.\\
{\bf Case 1.} $n$ is even.

In this case, we claim that $e_i'$ is not removed for all $1\leq i\leq n-2$. Suppose, to the contrary, that there exists an integer $i$ such that $e_i'$ is removed. In particular, we may assume that $i=1$.

Let $k=1,2,\dots, \frac{n-2}{2}-1$. We prove by induction on $k$ that for every pair of $(\Delta_{2k},\Delta_{2k+1})$, $e_{2k+1}$ is the unique edge in $E(\Delta_{2k})\cup E(\Delta_{2k+1})$ that is removed.

If $k=1$, since $e_1'$ is removed,  $e_2\in E(G[A])$. Otherwise, $G'=G+e_1'$ is a $K_4$-free $1$-planar graph, and $\lambda(G')>\lambda(G)$, a contradiction. If $e_2'$ is removed, then $G'=G+e_1'+e_2'-e_2$ is a $K_4$-free $1$-planar graph (refer to the left side of Figure \ref{trans1}).
\begin{figure}[h!]
    \centering
    \includegraphics[width=0.3\textwidth]{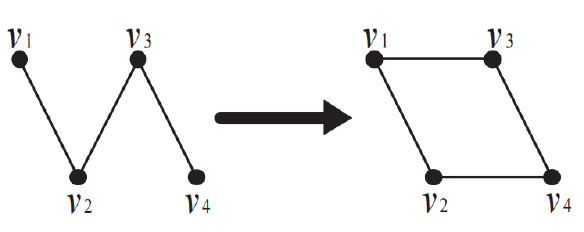}~~~~~~~~~~~~~\includegraphics[width=0.4\textwidth]{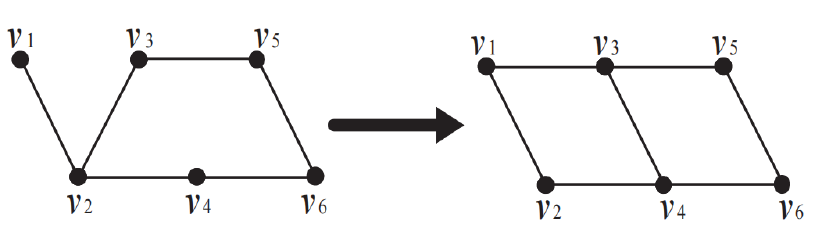}
    \caption{The transformation from $G$ to $ G'$.}
    \label{trans1}
\end{figure}
Moreover, by Claim \ref{cla-3} and $\lambda(G)\geq\sqrt{2n-4}$, we have
$$\mathbf{x}^{{\intercal}}A_{G'}\mathbf{x}-\mathbf{x}^{{\intercal}}A_{G}\mathbf{x}\geq 2\left( 2\times\frac{4}{(\lambda(G)-7)^2} -\frac{4}{\lambda(G)^2}\right)>0,$$
and thus $\lambda(G')>\lambda(G)$, a contradiction. Thus, $e_2'\in E(G[A])$ and $e_3$ is removed. Similarly, we have $e_3'\in  E(G[A])$. If $e_4$ is removed, $G'=G+e_3+e_1'-e_2$ is a $K_4$-free $1$-planar graph (refer to the right side of Figure \ref{trans1}), and $\lambda(G')>\lambda(G)$, a contradiction.

Now we can suppose that $k>1$ and that for any $j$ with $1\leq j < k$, $e_{2j+1}$ is the unique edge in $E(\Delta_{2j})\cup E(\Delta_{2j+1})$ that is removed. By induction, $e_{2k}\in E(G[A])$. Now we define a vector $\mathbf{y}$ as follow:
$$\mathbf{y}_x=\mathbf{y}_w=1,$$
\[
\mathbf{y}_{v_i} =
\begin{cases}
\mathbf{x}_{v_{2k+2}}, & i = 2, \\
\mathbf{x}_{v_{i-2}}, & i=4,6,\dots,2k+2,\\
\mathbf{x}_{v_i}, & \text{otherwise.}
\end{cases}
\]
Clearly, $\|\mathbf{y}\|_2=\|\mathbf{x}\|_2$. Now, if $e_{2k}'$ is removed, then let
$$G'=G+\sum_{i=1}^{k-1}e_{2i+1}-\sum_{i=1}^{k}e_{2i}+e_1'+e_{2k}'.~(\text{Refer to Figure \ref{trans2}}.)$$
\begin{figure}[h!]
    \centering
    \includegraphics[width=0.6\textwidth]{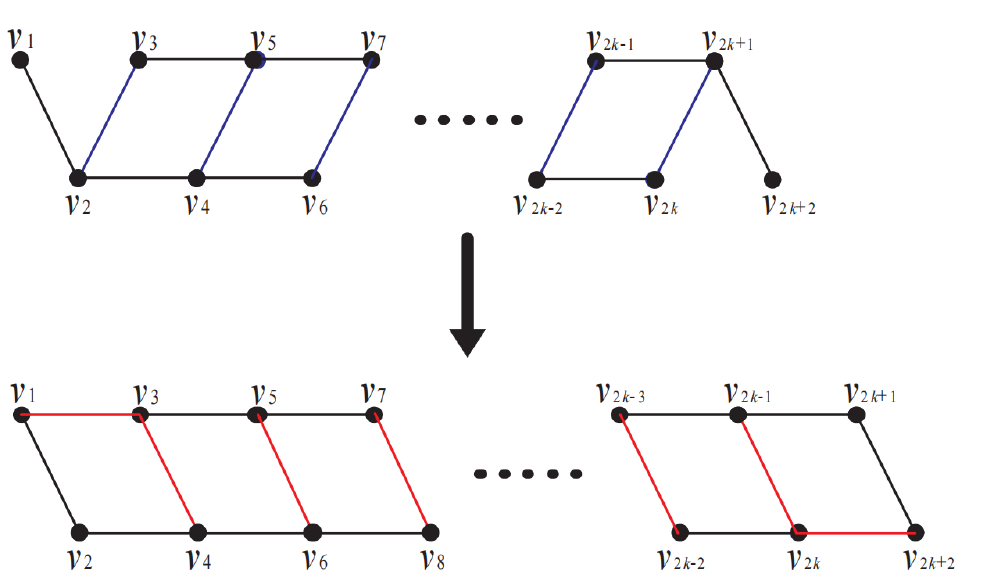}
    \caption{The transformation from $G$ to $ G'$. The blue edges indicate the edges removed during the transformation, while the red edges represent the edges added.}
    \label{trans2}
\end{figure}
It is easy to check that $G'$ is a $K_4$-free $1$-planar graph. Moreover, we have
\begin{align*}
\frac{1}{2}\left(\mathbf{y}^{{\intercal}}A_{G'}\mathbf{y}-\mathbf{y}^{{\intercal}}A_{G}\mathbf{y}\right)&=\sum_{i=1}^{k-1}(\mathbf{y}_{2i+1}\mathbf{y}_{2i+2})-\sum_{i=1}^k(\mathbf{y}_{2i}\mathbf{y}_{2i+1})+ \mathbf{y}_1\mathbf{y}_3+\mathbf{y}_{2k}\mathbf{y}_{2k+2}\\
&=\sum_{i=1}^{k-1}(\mathbf{y}_{2i+1}(\mathbf{y}_{2i+2}-\mathbf{y}_{2i}))-\mathbf{y}_{2k}\mathbf{y}_{2k+1}+ \mathbf{y}_1\mathbf{y}_3+\mathbf{y}_{2k}\mathbf{y}_{2k+2}\\
&\geq \sum_{i=1}^{k-1}(\mathbf{y}_{2i+1}(\mathbf{y}_{2i+2}-\mathbf{y}_{2i}))+\frac{8}{(\lambda(G)-7)^2}-\frac{4}{\lambda(G)^2}\\
&\geq \sum_{i=2}^{k-1}(\mathbf{y}_{2i+1}(\mathbf{y}_{2i+2}-\mathbf{y}_{2i}))+\frac{12}{(\lambda(G)-7)^2}-\frac{8}{\lambda(G)^2},
\end{align*}
and
\begin{align*}
\frac{1}{2}\left(\mathbf{y}^{{\intercal}}A_{G}\mathbf{y}-\mathbf{x}^{{\intercal}}A_{G}\mathbf{x}\right)&=\sum_{i=2}^{k+1}(\mathbf{y}_{2i}\mathbf{y}_{2i+2}-\mathbf{x}_{2i}\mathbf{x}_{2i+2})+\sum_{i=1}^{k+1}(\mathbf{y}_{2i}\mathbf{y}_{2i+1}-\mathbf{x}_{2i}\mathbf{x}_{2i+1})+\mathbf{y}_{2}\mathbf{y}_{n-2}+\mathbf{y}_{1}\mathbf{y}_{2}+\mathbf{y}_{2k+2}\mathbf{y}_{2k+3}\\
&-(\mathbf{x}_{2}\mathbf{x}_{n-2}+\mathbf{x}_{1}\mathbf{x}_{2}+\mathbf{x}_{2k+2}\mathbf{x}_{2k+3})\\
&\geq \sum_{i=2}^{k+1}(\mathbf{y}_{2i}\mathbf{y}_{2i+2}-\mathbf{x}_{2i}\mathbf{x}_{2i+2})+\sum_{i=1}^{k+1}(\mathbf{y}_{2i}\mathbf{y}_{2i+1}-\mathbf{x}_{2i}\mathbf{x}_{2i+1})+\frac{12}{(\lambda(G)-7)^2}-\frac{12}{\lambda(G)^2} \\
&= \mathbf{y}_{2}\mathbf{y}_{4}-\mathbf{x}_{2}\mathbf{x}_{4}+\mathbf{y}_{2k+2}\mathbf{y}_{2k+4}-\mathbf{x}_{2k+2}\mathbf{x}_{2k+4}
+\sum_{i=1}^{k+1}(\mathbf{y}_{2i}\mathbf{y}_{2i+1}-\mathbf{x}_{2i}\mathbf{x}_{2i+1})+\frac{12}{(\lambda(G)-7)^2}-\frac{12}{\lambda(G)^2}\\
&\geq \sum_{i=1}^{k+1}(\mathbf{y}_{2i}\mathbf{y}_{2i+1}-\mathbf{x}_{2i}\mathbf{x}_{2i+1})+\frac{20}{(\lambda(G)-7)^2}-\frac{20}{\lambda(G)^2}\\
&\geq \sum_{i=2}^{k-1}(\mathbf{y}_{2i+1}(\mathbf{y}_{2i}-\mathbf{x}_{2i}))+\frac{32}{(\lambda(G)-7)^2}-\frac{32}{\lambda(G)^2}\\
&\geq \sum_{i=2}^{k-1}(\mathbf{y}_{2i+1}(\mathbf{y}_{2i}-\mathbf{y}_{2i+2}))+\frac{32}{(\lambda(G)-7)^2}-\frac{32}{\lambda(G)^2}.
\end{align*}
Thus,
\begin{align*}
\frac{1}{2}\left(\mathbf{y}^{{\intercal}}A_{G'}\mathbf{y}-\mathbf{x}^{{\intercal}}A_{G}\mathbf{x}\right)&=\frac{1}{2}(\mathbf{y}^{{\intercal}}A_{G'}\mathbf{y}-\mathbf{y}^{{\intercal}}A_{G}\mathbf{y})+\frac{1}{2}(\mathbf{y}^{{\intercal}}A_{G}\mathbf{y}-\mathbf{x}^{{\intercal}}A_{G}\mathbf{x})\\
&\geq \frac{44}{(\lambda(G)-7)^2}-\frac{40}{\lambda(G)^2}\\
&>0,
\end{align*}

which implies that $\mathbf{y}^{{\intercal}}A_{G'}\mathbf{y}>\mathbf{x}^{{\intercal}}A_{G}\mathbf{x}$. Thus, $\lambda(G')>\lambda(G)$, a contradiction. Now we have $e_{2k}'\in E(G[A])$. Since at least one edge in $E(\Delta_{2k})$ is removed, we have that $e_{2k+1}$ is removed. Moreover, $e_{2k+1}'\in E(G[A])$. Therefore, it remains to show that $e_{2k+2}\in E(G[A])$. Suppose to the contrary and define
$$G'=G+\sum_{i=1}^{k}e_{2i+1}-\sum_{i=1}^{k}e_{2i}+e_1'.~(\text{See Figure \ref{trans3}.})$$
\begin{figure}[h!]
    \centering
    \includegraphics[width=0.6\textwidth]{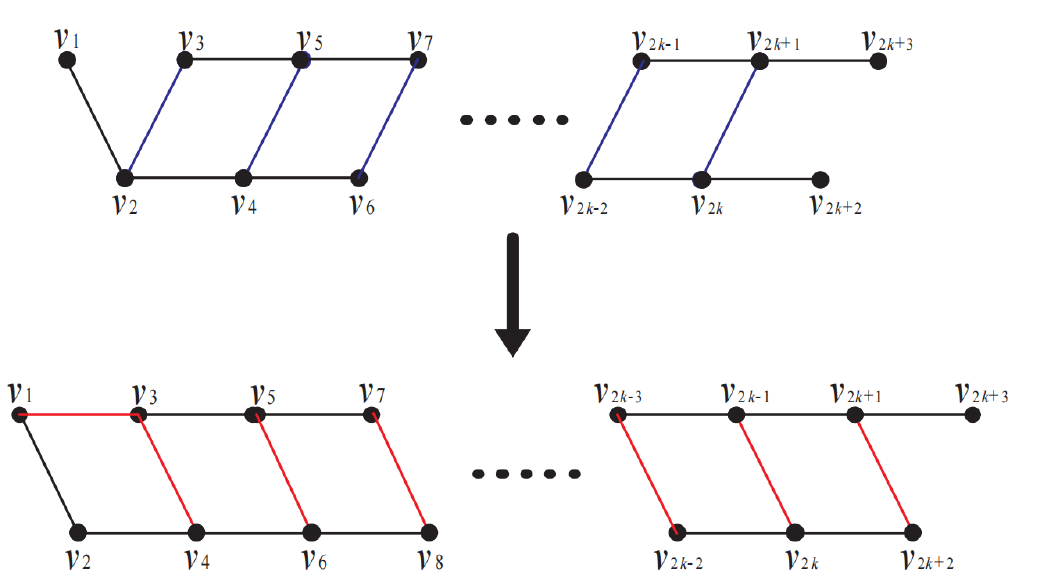}
    \caption{The transformation from $G$ to $ G'$. The blue edges indicate the edges removed during the transformation, while the red edges represent the edges added.}
    \label{trans3}
\end{figure}
It is easy to check that $G'$ is a $K_4$-free $1$-planar graph and similarly $\lambda(G')>\lambda(G)$, a contradiction. Therefore, $e_{2k+1}$ is the unique edge in $E(\Delta_{2k})\cup E(\Delta_{2k+1})$ that is removed,  and our induction is completed.

Next, we focus on $\Delta_{2n-2}$. Note that $E(\Delta_{n-2})=\{e_{n-2},e_{1},e_{n-2}'\}$. By induction, we have $e_{n-2}\in E(G[A])$.
If $e_1$ is removed, then simply adding the edge $e_1'$ to $G$, we can derive a contradiction in a manner similar to the previous one.
Thus, $e_{n-2}'$ is removed. Now we define $G'=G+e_1'+e_{n-2}'-e_1$. Clearly, $G'$ is $K_4$-free and $1$-planar. Moreover, $\lambda(G')>\lambda(G)$.
Now we have $e_{i}'\in E(G[A])$ for any $1\leq i\leq n-2$. If there exists an integer $1\leq i\leq n-2$ such that more than one edge in $E(\Delta_i)$ is removed.
Suppose without loss of generality that  $i=1$. By a proof that follows the same line as before, we can demonstrate a contradiction. Thus, $G[A]=C_{\frac{n-2}{2}}^{\square}$.\\
{\bf Case 2.} $n$ is odd.

 In this case, it is easy to check that there exists an integer $i$ such that $v_iv_{i+2}\notin E(G[A])~(i\mod n-2)$ since $G$ is $1$-planar. Thus, $G[A]$ is a subgraph of $C_{n-2}^{2-}$. Suppose that $v_{n-2}v_2\notin E(G)$, and define $\Delta_i$ to be the triangle in $C_{n-2}^{2-}$ formed by the three vertices $v_i$, $v_{i+1}$, and $v_{i+2}~(i\mod{n-2})$. We highlight that compared to Case 1, $\Delta_{n-2}$ cannot be defined.

The remaining part in this case is very similar to Case 1. For $\Delta_1$, if $e_1\notin E(G[A])$ or $e_1'\notin E(G[A])$, then we can use the same inductive method as in Case 1 to prove that for every pair of $(\Delta_{2k},\Delta_{2k+1})$ and $k=2,3,4,...,\frac{n-5}{2}$, $e_{2k+1}$ is the unique edge that is removed. Let $e^*$ and $e^{**}$ be the edges of $\Delta_1$ and $\Delta_{n-3}$ that is removed, respectively. Now, define
$$\mathbf{y}_x=\mathbf{y}_w=1,$$
\[
\mathbf{y}_{v_i} =
\begin{cases}
\mathbf{x}_{v_{2}}, & i = 2k+2, \\
\mathbf{x}_{v_{2i+2}}, & i=2,4,6,\dots,2k,\\
\mathbf{x}_{v_i}, & \text{otherwise.}
\end{cases}
\]
$$G'=G+\sum_{i=1}^{\frac{n-5}{2}}e_{2i+1}-\sum_{i=1}^{\frac{n-3}{2}}e_{2i}+e^*+e^{**}.$$
It is easy to check that $G'$ is a $K_4$-free $1$-planar graph, and $\mathbf{y}^{{\intercal}}A_{G'}\mathbf{y}>\mathbf{x}^{{\intercal}}A_{G}\mathbf{x}$. Thus, $\lambda(G')>\lambda(G)$, a contradiction.

Therefore, $e_2$ is the unique edge in $E(\Delta_1)$ that is removed. Then for $\Delta_2$, $e_2'\in E(G[A])$. Moreover, a derivation similar to that of $e_1$ shows that $e_3\in E(G[A])$.

Now we have for the pair of $(\Delta_{1},\Delta_{2})$, $e_{2}$ is the unique edge that is removed. Again, the same induction as in Case 1 shows that for every pair of $(\Delta_{2k-1},\Delta_{2k})$ and $k=1,2,3,...,\frac{n-3}{2}$, $e_{2k}$ is the unique edge that is removed. Therefore, $G[A]=C_{\frac{n-2}{2}}^{\square}$.
\endproof
\section{Proof of Theorem \ref{Th-K5}}\label{K5}
Let $G\in SPEX_{\mathcal{P}_1}(n,K_5)$ and $\mathbf{x}$ be the Perron vector of $G$ with  $\|\mathbf{x}\|_{\infty}=1$ and $x$ be a vertex such that $\mathbf{x}_x=1$. By Lemma \ref{lem-2}, there exist two vertices $x$ and $w$ such that $d_G(x)=d_G(w)\geq n-2$. Let $A=V(G)\setminus\{x,w\}$. Similarly to Claim \ref{cla-3} in Section \ref{K4}, we have $\frac{2}{\lambda(G)}\leq\mathbf{x}_v\leq \frac{2}{\lambda(G)-7}$ for all $v\in A$. Next, we consider two cases based on whether $xw\in E(G)$.

{\bf Case 1.} $xw \in E(G)$.

The proof of this case is very similar to Section \ref{K4}. Firstly, because $xw\in E(G)$, $G[A]$ is $K_3$-free. By a similar proof of Claim \ref{cla-2}, it is clear that there exists a $1$-planar drawing $D$ of $G$ such that $D(H)$ is a plane graph, where $H=K_{2,n-2}\subseteq G$. In $D(G)$, let the vertices in $A$ be labeled $v_1,v_2,\dots,v_{n-2}$ in counterclockwise order around $x$. Suppose without loss of generality that $xw$ is located in the region bounded by $xv_{n-2}wv_1$.
Since $D(G)$ is a $1$-plane graph and $D(H)$ is a plane graph, we have $N_G(v_i)\subseteq \{x,w,v_{i-2},v_{i-1}, v_{i+1}, v_{i+2}\}~(i\mod{n-2})$. Thus, $G[A]$ is a subgraph obtained from $P_{n-2}^{2+}$ by removing triangles.

The proof for the remaining part of this case differs from that in Section \ref{K4} in the following aspect: in Section \ref{K4}, regardless of whether $n$ is odd or even, the number of triangles that need to be removed is always even. However, in the current situation, when $n$ is odd, the number of triangles that need to be removed is odd, resulting in more than one possible extremal graph.

Firstly,  for $1\leq i\leq n-4$, we define $\Delta_i$ as the triangle in $P_{n-2}^{2+}$ formed by the three vertices $v_i$, $v_{i+1}$, and $v_{i+2}$. Let $e_i=v_{i}v_{i+1}$ and $e_i'=E(\Delta_i)-e_i-e_{i+1}=v_{i}v_{i+2}$. If $n$ is even, then we can follow the proof in Section \ref{K4} to show that for all $1\leq i \leq \frac{n-6}{2}$, $e_{2i+2}$ is the unique edge in $E(\Delta_{2i+1})\cup E(\Delta_{2i+2})$ that is removed. Therefore, $G[A]=QP_{n-2}$.

Now we suppose that $n$ is odd. If there exists an integer $1\leq i\leq n-4$ such that $e_i'$ is removed, then by the same induction as in  Section \ref{K4}, $i$ must be an odd integer and $e_{i-(2j+1)}$ is the unique edge in $E(\Delta_{i-(2j+1)})\cup E(\Delta_{i-(2j+2)})$ that is removed for all $0\leq j\leq \frac{i-1}{2}$. And $e_{i+(2j+2)}$ is the unique edge in $E(\Delta_{i+(2j+1)})\cup E(\Delta_{i+(2j+2)})$ that is removed for all $0\leq j\leq \frac{n-i-6}{2}$. (See Figure \ref{fig4} for example.) Thus, $G[A]\in \mathcal{P}_{n-2}^2$. Note that $e_1$ and $e_{n-3}$ correspond to $e_1'$ and $e_{n-4}'$, respectively, playing analogous roles in the above reasoning. Therefore, if  $e_1$ or $e_{n-3}$ is removed, then $G[A] \in \mathcal{P}_{n-2}^2 $ still holds.

\begin{figure}[h!]
    \centering
    \includegraphics[width=0.6\textwidth]{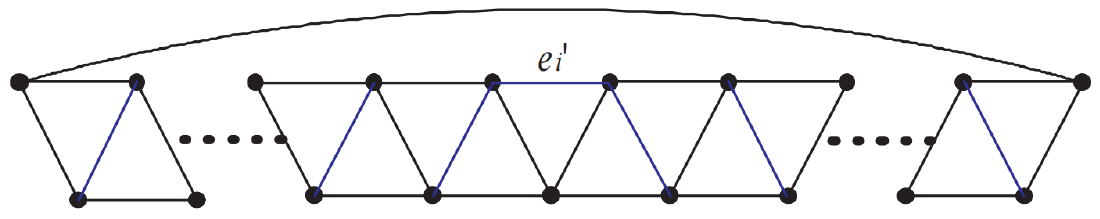}
    \caption{An example of $G[A]$. The blue edges indicate the edges removed.}
    \label{fig4}
\end{figure}
Now we assume that no edge in $\{e_1,e_{n-3},e_1',e_2',\dots,e_{n-4}'\}$ is removed. Let $E_1=\{e_2,e_4,e_6,...,e_{n-5}\}$ and $E_2=\{e_3,e_5,e_7,...,e_{n-4}\}$. Then there exists an integer $i$ such that in $E(\Delta_i)$, both $e_i$ and $e_{i+1}$ are removed. Again, by the same induction as in  Section \ref{K4}, $i$ must be an odd integer. Moreover, all edges removed are exactly those in the set $(\{e_j: j\leq i\}\cap E_1)\cup(\{e_j: j\geq i\}\cap E_2)$.  (See Figure \ref{fig5} for example.)  Therefore, we have $G[A]\in \mathcal{P}_{n-2}^2$. This completes the proof of Case 1.
\begin{figure}[h!]
    \centering
    \includegraphics[width=0.6\textwidth]{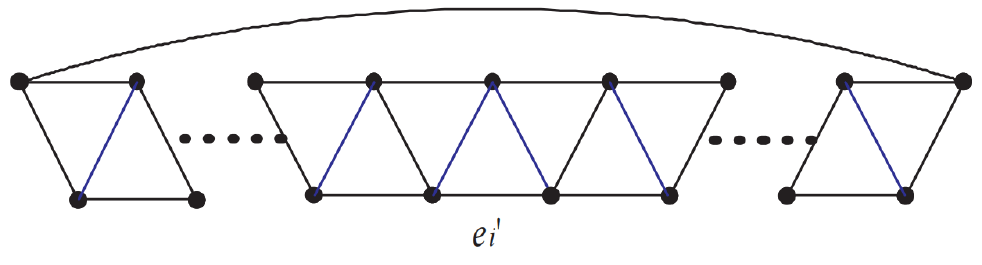}
    \caption{An example of $G[A]$. The blue edges indicate the edges removed.}
    \label{fig5}
\end{figure}

{\bf Case 2.} $xw \notin E(G)$.

Let $D$ be a $1$-planar drawing of $G$ that minimizes the crossing number and $H=K_{2,n-2}\subseteq G$. In $D(G)$, let the vertices in $A$ be labeled $v_1,v_2,\dots,v_{n-2}$ in counterclockwise order around $x$. If $xv$ crosses $wu$, then by Proposition \ref{pro-3}, no edge can cross $xu$, $wv$, and $uv$ (if $uv\in E(G)$). Thus, we have $xv_i$ crosses $wv_j$ only if $|i-j|=1$.

Suppose that there exist $v_i$ and $v_j$ such that $xv_i$ crosses $wv_j$. Then we have $|i-j|=1$. Assume without loss of generality that $i=1$ and $j=2$. Moreover, since $xv_1$ crosses $wv_2$, either $N_G(v_1)\cup N_G(v_2)\subseteq \{x, w, v_3,v_4\}$ or $N_G(v_1)\cup N_G(v_2)\subseteq \{x, w, v_{n-2},v_{n-3}\}$. We may suppose that $N_G(v_1)\cup N_G(v_2)\subseteq \{x, w, v_3, v_4\}$. Note that $G[A]$ is $K_4$-free, thus we have $e(G[\{v_1,v_2,v_3,v_4\}])\leq 5$ and $e(G[\{v_{n-2},v_{n-3},v_{n-4},v_{n-5}\}])\leq 5$. Now let $G'$ be the graph obtained from $G$ by removing all edges in $E(G[\{v_1,v_2,v_3,v_4\}])\cup E(G[\{v_{n-2},v_{n-3},v_{n-4},v_{n-5}\}])$, and adding edges $v_1v_2,v_1v_3,v_2v_3,v_2v_4,v_3v_4,v_{n-2}v_{n-3},v_{n-2}v_{n-4},v_{n-3}v_{n-4},v_{n-3}v_{n-5},v_{n-4}v_{n-5}$.We can modify $D$ in such a way that in $D(G')$, the subgraph $D(H[{v_{n-5},v_{n-4},v_{n-3},v_{n-2},v_1,v_2,v_3,v_4}])$ is a plane graph. (See Figure \ref{rd3}.) Define $G''=G'+v_{n-2}v_1$. Obviously, $G''$ is still $1$-planar and $K_5$-free. Moreover,
$$\mathbf{x}^{{\intercal}}A_{G''}\mathbf{x}-\mathbf{x}^{{\intercal}}A_{G}\mathbf{x}\geq 2\left( \frac{44}{(\lambda(G)-7)^2} -\frac{40}{\lambda(G)^2}\right)>0,$$
a contradiction.

Therefore, $D(H)$ is a plane graph. Thus, $G[A]=C_{n-2}^2$ if $n$ is even, and  $G[A]=C_{n-2}^{2-}$ if $n$ is odd, which completes our proof.
\section*{Acknowledgements}
This research was supported by the National Natural Science Foundation of China (Grant No. 12171089).

\end{document}